\documentclass[12pt,leqno]{article}
\usepackage[T1]{fontenc}
\usepackage{amsfonts}
\usepackage{amssymb}
\usepackage{mathrsfs}
\usepackage[latin1]{inputenc}
\usepackage{amsmath}
\usepackage[english]{babel}
 \usepackage{setspace}

\newcommand{\dis}{\displaystyle}

\newcommand{\lra}{\longrightarrow}

\newcommand{\N}{I\!\!N}
\newcommand{\R}{I\!\!R}

\newcommand{\fn}{\mbox{\boldmath$\it n$\unboldmath}}
\newcommand{\ssfn}{\mbox{\scriptsize\boldmath$\it n$\unboldmath}}

\newcommand{\fp}{\mbox{\boldmath$\it p$\unboldmath}}
\newcommand{\ssfp}{\mbox{\scriptsize\boldmath$\it p$\unboldmath}}

\newcommand{\fq}{\mbox{\boldmath$\it q$\unboldmath}}
\newcommand{\ssfq}{\mbox{\scriptsize\boldmath$ \it q$\unboldmath}}

\newcommand{\fnup}{{\fn}\,\mbox{-}\,\!\!\uparrow}
\newcommand{\fndown}{{\fn}\,\mbox{-}\,\!\!\downarrow}
\newcommand{\fnupdown}{{\fn}\,\mbox{-}\,\!\!\updownarrow}

\newcommand{\pcup}{\begin{array}[t]{c}\cup\\[-3.5ex]\cdot \end{array}}
\newcommand{\pbigcup}{\begin{array}[t]{c}\bigcup\\[-3.5ex]\cdot\end{array}}

\begin{document}

\thispagestyle{empty}

\begin{center}
 {\large \bf  On the compounding of higher order monotonic \\[1.5ex]
 pseudo-Boolean functions} \\[6ex]

  {\large \bf  Paul Ressel }
\end{center}

\enlargethispage{5.5cm}

\vspace*{1cm}

\textbf{\textsc{Abstract}} \\

  Compounding submodular monotone
(i.e.\ 2-alternating) set functions on a finite set preserves this property, as shown in 2010. A natural generalization to $k$-alternating functions was presented in 2018, however hardly
readable because of page long formulas. We give an easier proof of a more general result, exploiting known properties of higher order monotonic functions. \\



06E30 $\:   \cdot  \:$ 26A48 $\: \cdot  \:$ 26D07 $\:  \cdot  \:$ 26C99 \\[2ex]

{\large \bf 1. \  Introduction}  \\

Let $\: V \:$ be a finite non-empty set. A function $ \: \varphi : {\cal{P}}(V) \lra \R \:$ on $ \:{\cal P}(V),$ the set of all subsets of $ V ,$
i.e. a socalled \textit{pseudo-Boolean function}, is \textit{submodular}\/ if
\begin{equation} \label{1}
\varphi (A \cup \{ v \} ) - \varphi (A) \geq \varphi (B \cup \{ v \} ) - \varphi (B)
\end{equation}
for all $ \: A \subseteq B \:$ and all   $ v \in V \smallsetminus B \:.$ And $\:  \varphi  \:$
is \textit{increasing}\/  if $\:  \varphi (A) \leq \varphi (B)  \:$ whenever $ A \subseteq B \subseteq V \:.$
Condition \eqref{1}  has in many applications the interpretation that the marginal effect expressed by $\:  \varphi \:$ decreases for larger subsets
(the property of \glqq diminishing
returns\grqq). It is not surprising that submodular increasing functions are modelling many situations, both technical and social, for example
the influence in social networks. In this connection an interesting aggregation problem had been posed in \cite{Ke}: does \glqq local\grqq \  submodularity imply the corresponding
property  \glqq globally\grqq? This was confirmed 10 years later in \cite{Mo}. Now,
in another \glqq language\grqq, an increasing submodular function $\: \varphi \:$ on $\:{ \cal P}(V)
 \:$ is \glqq 2-alternating\grqq \ on $\:  \{ 0, 1\}  ^V \cong  {\cal P} (V)  \:,$ and it seems natural to consider the more general case of \glqq $k$-alternating\grqq \
 functions. For example, $\:   \varphi \:$ is 3-alternating
if, in addition to being increasing and submodular, the difference between the left and right hand side in
 \eqref{1} is further diminished
 if one more element is added. This idea is suggested in the recent work \cite{Ch}, whose central mathematical result (Theorem 4) however is given a very
 complicated and hardly readable proof, with page-long formulas. We shall give a much more transparent proof, based on existing theorems about higher order monotonic functions. Our result is also considerably more general.\\

  \vfill

\footnoterule

paul.ressel@ku.de

\newpage

 \textbf{Notations}.

 $ \N = \{ 1, 2, 3, \ldots \}, \N_0 = \{ 0, 1, 2, \ldots , \} , \R_+ = {[0, \infty[} ,  \, {\cal P}  (V) = $ set of all subsets
 of $\: V  \:,$  \newline
   $ 1_A (x) := \left\{\begin{array}{l}1, x \in A \\[-.8ex]
 0, x \not\in A \end{array} \right., \: [d] := \{ 1, 2, \ldots , d \}  \:$ for $\:  d  \in \N , \mathbf{1}_d  := (1, 1, \ldots , 1 )
  \in \N ^d \,,  $  \newline
  $  | \fn | := \sum^d_{ i = 1} n_i \:$ for  $\:  \fn \in \N^d_0,  \: | \alpha  | :=  $ cardinality of a finite set  $\: \alpha \: .$ \\
These two uses of the same symbol can hardly be mixed up; in fact, for $\:a = 1 _\alpha \in \{ 0, 1 \} ^d   \:$ we  have $\:  | a | = |  \alpha |  \:.$

  $ (f \times g ) (x, y ) := (f (x), g (y))  \:$ for mappings $\:  f, g, $ \\
   $  (f, g) (x) := (f( x), g (x)) \:$ for mappings $\:  f, g \:$ with the same domain\\[1ex]
  $
  \langle \sigma, \tau \rangle :=  \{ \gamma | \sigma \subseteq \gamma \subseteq \tau \} \:$ a \textit{set-interval}\/   where (usually) $\:  \sigma \subseteq \tau  \:$ \\[1ex]
 $  \langle  \sigma, \sigma \rangle = \{ \sigma \} $ is a special case \\[1ex]
  $\:  A  \pcup B, \! \pbigcup_{\!\!\!j}  \,  A_j $ for disjoint unions,  d.f. = distribution
  function (of some measure), \\
   $ x \odot y := (x_1 y_1, x_2 y_2, \ldots ) \:$ and $\: x  \vee  y : = (x_1 \vee y_1 , x_2 \vee y_2 , \ldots) $
for two vectors $\:  x, y  \:$ of equal dimension. \\[2ex]

{\large \bf 2. \ Multivariate higher order monotonicity} \\

Let $\:  I_1, \ldots , I_d \subseteq \R  \:$  be non-degenerate intervals, $\:  I := I _1 \times \ldots \times I_d  \:,$
and let $\:  f : I  \lra \R  \:$ be any function. For $\:  s \in I , h  \in \R^d_+  \:$
such that also $\:  s + h \in I  \:$ put
$$
(E_h f ) (s) := f (s + h )
$$
and $ \Delta_h := E_h - E_0\:,$ i.e. $ (\Delta_h \: f) (s) := f (s + h ) - f (s) \:,$
 and $\:  \nabla_h := - \Delta _h  \:.$ \\

  Since $\:  \{ E_h \}  \:$ is commutative (where defined), so are $\:  \{ \Delta _h \}  \:$
 and $\:  \{ \nabla _h \}  \:.$ In particular, with
 $\:  e_1, \ldots , e_d  \:$ denoting the standard unit vectors in $\:  \R^d , \Delta _{h_1 e_1}, \ldots,
 \Delta_{h_d e_d} \:$ commute. As usual $\:  \Delta ^0 _h \: f := f = : \nabla ^0_h \: f  \:.$ For
 $\:  \fn = (n_1, \ldots , n_d ) \in \N^d_0  \:$ and $\:  h \in \R^d_+  \:$
 we put
 $$
 \Delta ^{\ssfn}_h := \Delta ^{n_1}_{h_1 e_1} \Delta ^{n_2} _{ h_2 e_2 } \ldots \Delta ^{n_d} _{ h_d  e_d }
 $$
 and similarly $\: \nabla ^{\ssfn}_h \:.$ The multinomial theorem gives
 \begin{equation} \label{2}
 \left( \Delta ^{\ssfn}_h f \right) (s) = \sum_{ 0 \leq \mbox{\ssfq} \leq {\ssfn}} (-1)^{|{ \ssfn}| - |{ \ssfq}|}
{ \fn \choose \fq  }  f (s + \fq \odot h ) \: .
 \end{equation}
 Note that $\:  \Delta ^{\mathbf{1}_d}_h \neq \Delta_h  \:$ for $\:  d > 1  \:$  and $\: h \neq 0  \:.$ Also, $\:  \Delta ^{\ssfn} _h = 0 \:$
 if $\:  h_i = 0 < n_i \:$ for some $\:  i \leq d  \:.$ \\

\newpage

 {\bf Definition}.  $\:  f : I \lra \R \:$ is called
 \begin{itemize}
 \item[(i)] $\fn$-{\it increasing} (``$\fn\,\mbox{-}\!\!\uparrow$'') if
 $$
 ( \Delta ^{ \ssfp}_{ h} \:  f )  (s) \geq 0
 $$
 for all $\: s \in I , h \in \R ^d_+, \fp \in \N^d_0 , 0 \neq \fp \leq \fn \:$
 such that  $\:  s + \fp \odot h \in I  \:$

 \item[(ii)] $ \fn$-{\it decreasing} (``$\fndown$'') if instead
 $$
 \left( \nabla ^{ \ssfp}_h \:  f \right) (s ) \geq 0
 $$
   \item[(iii)] $\fn$-{\it alternating} (``$\fnupdown$'') if instead
  $$
  \left( \nabla^{ \ssfp}_h \:  f \right) (s) \leq 0 \: .
  $$
  \end{itemize}
  It is easy to see that,  using the notation $\:  (f ( - \cdot ) ) (s) := f (- s )  \:, $
 $$ \begin{array}{l} \dis
  f \: \mbox{ is }  \: \fndown \: \mbox{ on }  \: I \Longleftrightarrow f ( -  \cdot ) \: \mbox{ is }
 \: \fnup \: \mbox{ on } \: - I    \\[1.5ex]

 \dis f \: \mbox{ is } \: \fnupdown \: \mbox{ on } \: I \Longleftrightarrow - f (- \cdot ) \: \mbox{ is } \:
 \fnup \: \mbox{ on } \: - I \: .
 \end{array} $$

For $\:  \fn \in \{ 0,1 \} ^d  \:$ the $\:  I_j  \:$ considered here need not be intervals, just non-empty subsets of $\:  \R  \:, $ or even  $\:  \overline{\R} \:.$
Right-continuous bounded non-negative $\: \mathbf{1}_d\,\mbox{-}\,\!\!\uparrow \:$ functions on $\:  I  \:$ are precisely the distribution functions (\glqq d.f.s\grqq) of finite measures
on $\:  \overline{I}  \:$ (closure in $\:  \overline{\R}^d) \:,$ see \cite{Re1} Theorem 7, a  result which will be used later on. Functions which are
$\: \mathbf{1}_d\,\mbox{-}\,\!\!\uparrow (\downarrow, \updownarrow) \:$  are also called \textit{fully}\/ $d$-\textit{increasing} (-\textit{decreasing}, -\textit{alternating}), and this   notion will  now be extended: \\

\textbf{Definition}. \ Let $\:  I_1, \ldots , I_d \subseteq \overline{\R}  \:$ be any non-empty subsets, $\:  I := I_1 \times \ldots \times I_d, 1 \leq k \leq d \:.$
Then $\:  f :  I \lra \R  \:$ is \textit{fully} $ k$-\textit{increasing} \ (\glqq$\mathbf{1}_k\,\mbox{-}\,\!\!\uparrow$\grqq) \
 iff
  $$
  ( \Delta ^{\ssfp}_h \, f ) (s) \geq 0 \; \mbox{ for each  } \; 0 \lvertneqq \mathbf{p} \leq \mathbf{1}_d  \: \mbox{ with } \: | \mathbf{p} | \leq k
  $$
  and for each $\:  s \in I  \:$ and $\:  h \in \R^d _+ \:$  such that $ \: \{ s + \mathbf{q} \odot h \mid \mathbf{q }\in \{ 0, 1 \} ^d , \: \mathbf{q} \leq \mathbf{p} \} \subseteq I  \:. $ \\

If instead $\:  ( \nabla ^{\ssfp} _ h f ) (s) \geq 0  \:$ we call  $\:  f  \:$ \textit{fully} $k$-\textit{decreasing}  \ $(\glqq\mathbf{1}_k\,\mbox{-}\,\!\!\downarrow$\grqq) \,
and if $\:  ( \nabla ^{\ssfp}_h f ) (s) \leq 0  \:,$ $\:  f  \:$ is by definition \textit{fully} $\:  k $-\textit{alternating}
\ (\glqq $ \mathbf{1}_k\,\mbox{-}\,\!\!\updownarrow$\grqq) . \\

For the important special case where $\:  I_j = \{ 0, 1 \}  \;  \forall \: j \leq d  \:,$ i.e. for  (\glqq pseudo-Boolean\grqq) functions
on $\:  \{ 0, 1 \} ^d  \:, $ it is  sometimes useful to identify $\:  \{ 0, 1 \}^d  \:$ with $\:  {\cal P} ([d]) := \{ \alpha | \alpha \subseteq [d] \} \:.$
Since $\:  \Delta ^1_0 = 0  \:$ and $\: \Delta ^0 _1 = \mbox{id} = \Delta^0_0  \:,$ only $\:  \Delta^{\ssfp}_h  \:$ with
$\:  h = \fp \in \{ 0, 1 \} ^d \smallsetminus \{ 0 \}  \:$ have to be considered.
It is then reasonable to use the simplified notation
$$
\Delta _\alpha := \Delta ^a _a \: \mbox{ for  } \: a = 1_\alpha \in \{ 0, 1 \} ^d \smallsetminus \{ 0 \}
$$
(complemented by $\:  \Delta_\emptyset = \mbox{id}) \:.$ \\

We write likewise  $\:  E_\alpha : = E_a  \:.$ Both $\: \Delta_\alpha f  \:$ and $\:  E_\alpha f  \:$ have the domain
$\:  \{ \gamma \subseteq [d] | \gamma \subseteq \alpha ^c \} \:,$ and for $\:  \gamma \subseteq \alpha ^c  \:$
$$ ( \Delta _ \alpha f ) (\gamma) = (( E_\gamma \Delta _\alpha ) (f)) ( \emptyset) = (( \Delta _\alpha E_\gamma ) (f)) ( \emptyset).
$$
Clearly $\: \Delta _\alpha \circ \Delta _\beta = \Delta _{\alpha \cup \beta} \:$ for disjoint $\: \alpha , \beta  \:.$
Note that
$$
(\Delta _\alpha f ) (\emptyset ) = f ( \alpha ) - \sum_{{\gamma \subseteq \alpha \atop | \gamma | = | \alpha | - 1} } f ( \gamma)  + \sum_{ { \gamma \subseteq \alpha \atop
| \gamma | = | \alpha | - 2 }} f (\gamma ) \mp \ldots + (-1) ^{ | \alpha  |} f ( \emptyset)  \:.
$$
The following identity (for $\:  x _1, \ldots , x_d \in \R )  \:$
$$
\prod ^d _{ i = 1 } x_i = \prod^d _{i = 1} [(x_i - 1 ) + 1] = \sum _{ \alpha \subseteq [d]} \prod _{ i \in \alpha} (x_i - 1) \qquad
\left( \prod_\emptyset := 1  \right)
$$
holds of course also within the commutative algebra generated by $ \: \{ E_{\{ i \} } | i \in [d]\} ,\:$  and leads to
$$
\sum_{ \alpha \subseteq [d]} \Delta _\alpha = \prod ^d _{ i = 1} E_{\{ i \}} = E_{[d]} \:,
$$
i.e. to $\:  \sum _{\alpha \subseteq [d]} ( \Delta _\alpha f ) (\emptyset ) = f ( [d]) \:.$ \\

Slightly more general, and of importance later on, for $\:  \beta \subseteq \gamma \subseteq [d] \:$
\begin{equation} \label{3}
\sum_{ \alpha \in \langle \beta, \gamma \rangle} \Delta _\alpha = \Delta _\beta \sum_{ \alpha \subseteq \gamma \smallsetminus \beta} \Delta_\alpha
= \Delta _\beta E_{ \gamma \smallsetminus \beta } \: .
\end{equation}

We mention that fully $k$-alternating pseudo-Boolean functions are called  \ \glqq$AD-k$\grqq \  in \cite{Ch}. \\[2ex]

\newpage

 \textbf{\large 3.\ Multilinear polynomials} \\

 Any (pseudo-Boolean) function  $\:  f : \{ 0, 1 \}^d \lra \R  \:$
  has an extension $\:  \tilde{f} \:$ to a socalled \textit{multilinear polynomial}
  \begin{equation}
  \label{4}
   \tilde{f} (x) := \sum_{\alpha  \subseteq [d]} f (\alpha )  x^\alpha (\mathbf{1} - x )^{\alpha^c}  \:, \quad x \in \R^d
   \end{equation}
   where we use the abbreviations $\: x ^\alpha := \prod_{i \in \alpha } x_i , x ^\emptyset := 1  \:$ and $\:  \mathbf{1} := \mathbf{1}_d \:.$
   \glqq Multilinear\grqq \ means here that no variable appears in a power $\:  > 1 \:$ in $\: \tilde{f}  \:;$ $\:  \tilde{f} \:$ is therefore an affine
   function of each variable $\:  x_i  \:.$ Note that $\:  \tilde{f} (1_\alpha) = f (\alpha) \: \forall \; \alpha \subseteq [d]  \:,$
   so $\:  \tilde{f}   \:$ is uniquely determined; or in other words, each multilinear polynomial is the extension of its restriction to
   $\: \{ 0, 1 \} ^d  \:, $ where we freely identify $\: \alpha \subseteq [d] \:$ with $\:  1 _\alpha \in \{ 0, 1 \} ^d  \:.$
   It is immediate that $\:  f \geq 0  \:$ iff $\: \tilde{f} \mid [0, 1] ^d \geq 0  \:.$ \\

   Let for $\:  \emptyset \neq \beta \subseteq [d]  \:$ the partial derivative of $\:  \tilde{f} \:$ w.r. to $\:  x_i, i \in \beta  \:$ be $\: \partial ^\beta
   \tilde{f} \:$ (every other partial derivative of $\:  \tilde{f} \:$  is obviously 0). Then for any $\:  p  \in [d] \:$
    $$
    ( \partial ^{\{ p \} } \tilde{f})  (x) = \sum_{ \alpha \subseteq [d]\smallsetminus \{ p \} } ( \Delta _{ \{ p \}} f) (\alpha)
    x^\alpha (\mathbf{1} - x ) ^{[d] \smallsetminus ( \alpha \cup \{ p \} )}
    $$
    by an application of the product role, i.e. $\:  \partial ^{\{ p \}}  \tilde{f}\:$  is multilinear in $\:   x_i  \:$ for $\:  i  \in [d]  \smallsetminus \{ p \} \:.$
    By iteration we obtain for any $\:  \emptyset \neq \beta \subseteq [d] \:$
    \begin{equation} \label{5}
    ( \partial ^\beta  \tilde{f} ) (x) = \sum _{ \alpha \subseteq [d] \smallsetminus \beta } ( \Delta _\beta  f ) ( \alpha )
    x ^\alpha ( \mathbf{1} - x ) ^{[d]  \smallsetminus ( \alpha \cup \beta )}
    \end{equation}
    including finally
    $$
    ( \partial ^{[d]} \tilde{f} ) (x) = ( \Delta _{[d]} f ) ( \emptyset) , \quad \mbox{ a constant.}
    $$
    That is, $\:  \partial ^\beta \tilde{f} \:$ is the multilinear extension of $\:  \Delta _\beta f  \:$ on $\:  \{ 0, 1 \} ^{\beta ^c} \:.$\\

    Now  \eqref{5}  implies
      \begin{equation} \label{6}
       (\partial ^\beta \tilde{f} ) (0) = ( \Delta _ \beta f) ( \emptyset) , \quad  \beta \subseteq [d] ,
       \end{equation}
       and in the likewise \glqq canonical\grqq \ representation
       $$
       \tilde{f}  (x) = \sum_{ \alpha \subseteq [d]} c_\alpha x^\alpha $$
       we have obviously $\:  c_\alpha =  ( \partial ^\alpha \tilde{f} )(0) \:.$ Combining this with \eqref{6}
       we get
       \begin{equation} \label{7}
       \tilde{f} (x)  = \sum_{\alpha} ( \Delta _\alpha f ) ( \emptyset) x^\alpha  .
       \end{equation}
       We'll need later on the following result: \\

   \textbf{Lemma 1}. \ \textit{For } $\:  f : \{ 0 , 1 \} ^d \lra \R  \:$ \textit{and its multilinear extension} $\:  \tilde{f} \:$ \textit{we have}
   $$
   f  \;  is \; \mathbf{1}_k \,\mbox{-}\!\,\!\uparrow  \:  \Longleftrightarrow  \: \tilde{f} \; is \; \mathbf{1}_k\,\mbox{-}\,\!\!\uparrow \: on \: [0,1]^d.
   $$

   \textbf{Proof}. \ $ \tilde{f} \:$ is a polynomial, in particular $\:  C^\infty \:.$ Therefore $\:  \tilde{f}  \:$
   is $\:  \mathbf{1}_k\,\mbox{-}\,\!\!\!\uparrow \:$ (on $\:  [0,1]^d ) $ if and only if
   $$
   (\partial^ \beta \tilde{f} )(x) \geq 0  \quad \forall \; | \beta | \leq k , \; \forall \: x \in [0,1]^d
   $$
   which, as we just saw, is equivalent with
   $$
    ( \Delta _\beta f ) ( \alpha ) \geq 0 \quad \forall \: | \beta | \leq k, \; \forall \: \alpha \subseteq \beta ^c  \: ,
    $$
    the defining property of $\:  f  \:$ being $\: \mathbf{1}_k\,\mbox{-}\,\!\!\uparrow  \:. \  \Box$ \\

    \textbf{Example 1}. \ For $\:  d =
      3, k = 2  \:$ consider $\:  f : {\cal P} ( [3]) \lra \R  \:$ given by $\:  f (\alpha) :=
    | \alpha | \vee 1  \:.$ Then
    $$
    \tilde{f} (x) = 1 + x _1 x_2 + x_1  x_3 + x_2 x_3 - x_1 x_2 x_3
    $$
    and
    $$ \begin{array}{l}
    ( \partial ^{\{1 \} } \tilde{f} ) (x) = x_2 + x_3  - x_2 x_3 \; \mbox{ etc.} \\[1.5ex]
   (\partial ^{\{ 1, 2 \} } \tilde{f} ) (x) = 1 - x_3  \; \mbox{ etc.}
   \end{array}
   $$
   are all non-negative on $\:  [0,1]^3 \:; $ however
    $$ \partial ^{ \{ 1, 2, 3 \} }  \tilde{f} = - 1 \:,
    $$
   showing $\:  f  \:$ to be $\: \mathbf{1}_2 \,\mbox{-}\,\!\!\uparrow \:,$ but not $\:  \mathbf{1}_3\,\mbox{-}\,\!\!\uparrow\:.$
   Slightly more general, $\:  f ( \alpha ) := | \alpha | \vee 1 \:$ is for any $\:  d \geq 3 \; \: \mathbf{1}_2\,\mbox{-}\,\!\!\uparrow \:$
   and not $\:  \mathbf{1}_3\,\mbox{-}\,\!\!\uparrow \::$
   we have
   $$
   \tilde{f} (x) = \sum^d _{ i = 1} x_i + \prod^d _{ i = 1} ( 1 - x_i ) \: ,
   $$
   whence
   $$ \begin{array}{l} \dis
   ( \partial ^{\{ i \}} \tilde{f}) (x) = 1 - \prod_{ \ell \neq i}     ( 1 - x _\ell ) , \\  [1.5ex]

   \dis  (\partial ^{\{ i, j \}}  \tilde{f} ) (x) =  \prod_{ \ell \neq i, j} (1 - x_\ell ) \quad \mbox{ for } \; i \neq j
 \: ,
\end{array} $$
and
$$
 \partial ^\alpha \tilde{f} \equiv  - 1 \quad \mbox{ for }  \: | \alpha | = 3 \:  .
 $$

 \vspace*{.4cm}

 \textbf{\large 4. \ A combinatorial intermezzo} \\

The following Lemma (of combinatorial nature) will play a crucial role in the proof of the main result. We shall use \textit{set-intervals} in
 $\:{\cal P}([d]) \:$ of the form
 $$
 \langle \sigma, \tau \rangle := \{ \gamma \mid \sigma \subseteq \gamma \subseteq \tau \}  \:,
 $$
 including as a special case singletons
 $$
 \langle \sigma \rangle := \langle \sigma , \sigma \rangle  = \{ \sigma \} \: .
  $$
  Note that
 $\: \langle \sigma, \tau \rangle \neq \emptyset \:$ iff $\: \sigma \subseteq \tau \:,$ and that
 $$
 \langle \sigma_1, \tau_1 \rangle \cap \langle \sigma_2, \tau_2  \rangle = \langle  \sigma _1 \cup \sigma_2 , \tau_1 \cap \tau_2 \rangle \: .
 $$

 \textbf{\large Lemma 2}. \ \textit{Let}  $\:  k, d \in \N , k  \leq d  \:,$  \textit{and} $\:  x_1, \ldots, x_d \in \R^k  \:.$
 \textit{For non-empty} $\:  \alpha, \beta \subseteq [d] \:$ \textit{define}
 $$
 \alpha \sim \beta : \Longleftrightarrow \max_{ i \in \alpha} x_i = \max_{ i \in \beta} x_i  \quad ( \in \R^k ).
 $$
 \textit{Then} $\:   \{ \gamma \subseteq [d] \mid   \, | \gamma | \geq k \}  \:$  \textit{is the disjoint union of set-intervals}
 $\:  \langle \sigma_j , \tau_j \rangle  \:$ \textit{with} $\:  | \sigma_j | = k , \sigma_j \subseteq \tau_j  \:$  \textit{and}
 $\:  \sigma_j \sim \tau_j  \:$ \textit{for each} $\:  j  \:.$ \\

 \textbf{\large Proof}. \ For $\:  k = 1 \:$ we may assume $\:  x_1 \leq x_2 \leq \ldots \leq x_d  \:,$ and then
 $$
 \{ \gamma \subseteq [d] \mid \, | \gamma | \geq 1 \} = \langle \{ d \} , [d] \rangle \pcup  \langle \{ d - 1 \} , [d - 1] \rangle
 \pcup \ldots \pcup \langle \{ 2 \}, [2] \rangle   \pcup  \langle \{ 1 \} \rangle
 $$
 has the required properties.\\

 For $ \: k \geq 2 \:$ and  $\:  d = k + 1  \:$  choose $\:   \alpha \subseteq [d]  \:$  of size $\:  k  \:$ such that
 $\:   \alpha \sim [d] \:$ (which is evidently possible for any  $\:     d > k  \:).$ Then
 $$
 \{ \gamma \subseteq [d] \mid | \gamma | \geq k \} = \langle \alpha, [d] \rangle \pcup \!\! \pbigcup_{\!\!\!a \in \alpha} \:  \langle [d] \smallsetminus \{ a \} \rangle \: .
 $$
 We now proceed by induction and suppose  the result to be true for some $\:  d \geq 3  \:$ and each $\:  k \leq d  \:.$ Let $\:  x_ 1, \ldots , x_{d + 1}  \in \R ^k \:$
 be given, $\: k \geq 2 \:.$ It is no restriction to assume $\:  k < d  \:$ and
 $$
 x_{ d + 1} (k) = \max_{i \leq d + 1} x_i (k)\:.
 $$
 By assumption,
 $$
 \{ \gamma \subseteq [d] \mid | \gamma | \geq k \} = \!\!\pbigcup_{\!\!\!j}  \; \langle \xi _j, \eta _j \rangle
 $$
 is the disjoint union of set-intervals, with $\: \xi _j \subseteq \eta_j \subseteq [d] , | \xi_j | = k  \:$ and
 $\:  \xi_j \sim \eta_j  \:$ for each $\:  j  \:.$ Let $\:  y_i := (x_i (1), \ldots , x_i (k - 1 )) \in \R ^{ k - 1}  \:$
 be the projection of $\:  x_i , i = 1, \ldots , d  \:.$ Making use once more of the induction hypothesis we have
 $$
  \{ \gamma \subseteq [d]  \mid \, | \gamma | \geq k - 1 \} = \!\!\pbigcup_{\!\!\!p}  \: \langle \alpha_p, \beta_p \rangle
  $$
  with $\:  \alpha _p \subseteq \beta_p  \subseteq [d] , | \alpha_p | = k - 1 , \alpha_p \approx \beta_p  \:$ for all $\: p,  \:$
  where $\:  \alpha \approx \beta \:$ means $\: \max\limits _{ i \in \alpha} y_i = \max\limits_{ i \in \beta } y_i  ( \in \R^{k - 1} ). \:$ \\

  We now put
  $$
  \overline{\alpha}_p := \alpha_p \cup \{ d + 1 \} , \quad \overline{\beta}_p := \beta _p \cup \{ d + 1 \} \: ,
  $$
  then  $\:  | \overline{\alpha}_p | = k , \, \overline{\alpha}_p \subseteq \overline{\beta}_p \subseteq [d + 1] \:,$ and $\:  \overline{\alpha}_p \sim \overline{\beta}_p \:,$
  since  for $\:  \ell <  k  \:$
  $$  \begin{array}{lcl}
 \dis  \max_{ i \in \overline{\alpha}_p} x_i ( \ell)    & =  &  \dis \left( \max_{ i \in \alpha_p} x_i ( \ell)\right) \vee x_{d + 1} (\ell)  \\[1.5ex]

 & = & \dis \left( \max_{ i \in \alpha_p} y_i (\ell)\right) \vee x_{d + 1} (\ell)  \\[1.5ex]

 & = & \dis  \left( \max_{ i \in \beta_p} y_i ( \ell) \right) \vee x_{d + 1} (\ell)  \\[1.5ex]

 & = & \dis \max_{i \in \overline{\beta}_p} x_i ( \ell)
\end{array}  $$
and
$$
\max_{i \in  \overline{\alpha}_p} x_i ( k ) = x_{d + 1} (k) = \max_{ i \in \overline{\beta}_p} x_i ( k) \: .
$$
For any $\: j \:$ and $\:  p  \:$ we have
$$
 \langle \xi _j, \eta_j \rangle  \cap \langle \overline{\alpha}_p, \overline{\beta}_p \rangle =
 \langle \xi_j \cup \overline{\alpha}_p, \eta_j \cap \overline{\beta}_p \rangle = \emptyset
 $$
 since $\:  d + 1 \in \overline{\alpha}_p \:,$ but $\:  d + 1  \not\in \eta_j  \:.$ \\

 For $\:  p \neq q   \:$ likewise
 $$
 \langle \overline{\alpha}_p, \overline{\beta}_p \rangle \cap \langle \overline{\alpha}_q, \overline{\beta}_q \rangle = \emptyset
 $$
 because otherwise $\:  \alpha _p \cup \alpha_q \subseteq \beta_p \cap \beta_q  \:,$ contradicting the choice of $\:  \alpha_p , \beta_p  \:.$ \\

 So, finally
 $$
 \{ \gamma \subseteq [d + 1]| \, | \gamma | \geq k \} = \!\!\pbigcup_{\!\!\!j} \langle \xi_j , \eta_j \rangle \pcup \!\!\pbigcup_{\!\!\!p} \langle \overline{\alpha}_p,
 \overline{\beta}_p  \rangle
 $$
 is a partition into disjoint set-intervals as claimed. $ \Box$ \\[2ex]

 \textbf{\large 5. \  The main result}. \\

 In  \cite{Ch}, Theorem 4 the following is shown: let $\:  f : \{ 0, 1 \} ^d \lra [0,1] \:$ and $\:  g_1, \ldots , g_d : \{ 0, 1 \} ^k \lra [0,1]  \:$
 be all fully $k$-alternating   (\glqq$\mathbf{1}_k\,\mbox{-}\!\!\,\updownarrow$\grqq), then also their \glqq compounding\grqq \
 $\:  h : \{ 0, 1 \} ^k \lra \R \:,$ defined by
 $$
 h (x) := \sum_{ \alpha \subseteq [d]} f (\alpha ) \prod _{ i \in  \alpha} g_i ( x) \prod  _{ j \in \alpha ^c } (1 - g _j (x))
 $$
 has this property. The proof there is based on the multilinear extensions of $\:  f , \{ g _i \}  \:$ and $\:  h  \:,$ but it is hardly readable, with
 formulas longer than a page. Since the result is true (see below), I believe their proof is, too, although I didn't check it in detail --- by lack of patience.  \\

 We will prove a more general result, allowing $\: g_1, \ldots , g_d  \:$ to be any $\:  \mathbf{1}_k\,\mbox{-}\!\!\,\updownarrow \:$ functions on an arbitrary
 product subset of $\:  \overline{\R}^k  \:.$ Only $\: f \:$ has to remain a pseudo-Boolean function.  We  shall first deal with
 $\:  \mathbf{1}_k\,\mbox{-}\!\,\uparrow \:$ functions (generalizing increasing supermodular functions), and then deduce from it
 the statement about $\: \mathbf{1}_k\,\mbox{-}\!\!\,\updownarrow \:$ functions in a straightforward way.\\

 We shall need the following approximation result.\\

 \textbf{Lemma 3}. \ \textit{Let $\:  A = A_1 \times \ldots \times A_ k  \:$ be a product of non-empty subsets $\: A_j \subseteq \overline{\R}  \:,$
 and let $\:  g : A \lra [0,1] \:$ be $\:  \mathbf{1}_k\,\mbox{-}\!\!\uparrow \:$ and such that $\:  \sup g (A) = 1  \:.$
 Then there is a net $\:  (g_\alpha)  \:$ of distribution functions of probability measures with finite support contained in $\: A  \:,$
 which converges pointwise to $\:  g  \:.$ } \\

 {\bf Proof.} \ Let $\: \alpha _j \subseteq A_j  \:$ be finite and non-empty, $\:  1 \leq j \leq k  \:,$ and
 $\:  \alpha := \alpha _1 \times \ldots \times \alpha_k  \:;$ we may assume the $\:  \alpha _j  \:$ so large that
 $\:  g ( \max \alpha ) >  0  \:.$ The restriction $\:  g \mid \alpha   \:$ is $\: \mathbf{1}_k\,\mbox{-}\,\!\!\uparrow \:,$
 (automatically right-continuous on $\:  \alpha (!)  \:),$ and so there exists  by \cite{Re1}, Theorem 7 a finite measure $\:  \nu_\alpha  \:$ on $\:  \alpha  \:$ with
 d.f. $\:  g \mid \alpha \:.$ We have $\: \nu_\alpha ( \alpha ) = g ( \max \alpha ) >  0  \:,$ hence $\:  \mu_\alpha := \nu _\alpha / g ( \max \alpha )  \:$
 is a probability measure on $\:  \alpha  \:,$ which is extended trivially to a probability measure on $\: \overline{A}  \:,$ with
 $\:  \mu_ \alpha ( \overline{A} \smallsetminus \alpha ) := 0  \:.$ By $\: g _\alpha : A \lra [0,1]  \:$ we denote the d.f. of this extended
 $\:  \mu_\alpha \:.$\\

  In order to see that $\: g_\alpha    \:$ converges pointwise to $\:  g  \:, $ let $\:  0 < \varepsilon <   1/2 \:$ and some
 (finite, non-empty) product set $\:   \alpha _0 \subseteq A  \:$ be given. Choose $\:  \alpha \supseteq \alpha _0  \:$ (a product set, too)
 so large, such that $\:  g ( \max \alpha ) \geq 1 - \varepsilon \:.$ Then for any $\:  a \in \alpha  \:$
 $$ \begin{array}{lcl}
 \dis \left|  g_\alpha ( a) - g (a) \right| &  = & \dis  \left| \frac{g (a)}{g ( \max \alpha ) } - g ( a ) \right| = g (a) \cdot
 \frac{1 - g (\max \alpha )}{g (\max \alpha)} \\[2ex]

 & \leq & \dis g (a)  \frac{\varepsilon}{1 - \varepsilon} \leq 2 \varepsilon \:. \end{array}
 $$
 Noting that the family of finite product sets in $\:  A  \:$ is upwards filtering, the proof is complete. $  \ \Box $  \\

 \textbf{Theorem 1}. \ \textit{Let  $\:  k, d \in \N , k \leq d , \emptyset \neq A_j \subseteq \overline{\R} \:$ for
 $\:  j = 1, \ldots , k , A := A_1 \times \ldots \times A_ k  \:.$ Let $\: g_i : A \lra [0,1]  \:$ for $\:   i = 1, \ldots , d  \:$ and
 $\:  f : \{ 0, 1 \} ^d  \lra  \R  \:$ be given. Define $\:   h : A \lra \R  \:$ by }
 $$
 h (x) := \sum _{ \alpha \subseteq [d]} f (\alpha ) \prod_{ i \in \alpha} g_i (x) \prod_{ j \in \alpha ^c } [1 - g_j (x)]  \: .
 $$
 \textit{Then, if $\:  g_1, \ldots g_d  \:$ and $\:  f  \:$ are all $\:  \mathbf{1}_k\,\mbox{-}\,\!\!\uparrow \:,$  so is $\:  h  \:.$ } \\

 \textbf{Proof}. \ With $\:  \tilde{f} \:$ as the multilinear extension of $\:  f  \:,$ and $\:  g : = (g_1, \ldots , g_d ) : A \lra [0,1]^d  \:,$ we have
 $\:  h = \tilde{f} \circ g  \:,$ and by Lemma 1 $\:  \tilde{f} \:$ is also $\:  \mathbf{1}_k\,\mbox{-}\,\!\! \uparrow\:$ on $\:  [0,1]^d \:.$ \\

 We first consider the case that $\:  g_i  \:$ is the d.f. of some one-point measure $\:  \varepsilon_{a_i} \:,$ where
 $\:  a_i \in A  \:.$ Then
 $$
 g_i = 1 _{[a_i, \infty] \cap A}
 $$
 and for $\:  \emptyset \neq \alpha \subseteq [d] \:$
 $$
 \prod_{i \in \alpha} g_i = 1 _{[\max\limits_{i \in \alpha} a_i, \infty] \cap A} \:,
 $$
 and then by (7)
 $$ \begin{array} {lcl}
 \dis \tilde{f} \circ g & = & \dis \sum_{ \alpha \subseteq [d]} ( \Delta_\alpha f ) (\emptyset) \cdot \prod_{i \in \alpha} g_i  \\[2ex]

 & = & \dis \sum_{\alpha \subseteq [d]} ( \Delta  _\alpha  f ) ( \emptyset) \cdot 1 _{[\max_{i \in \alpha }  a_i , \infty ] \cap A  } \:.
 \end{array} $$

 For $\:  k = d  \:$  we have $\:   ( \Delta_ \alpha f ) ( \emptyset)    \geq 0  \:$ for each
 $\:   \alpha \subseteq [d] \:,$ implying directly that $\:  \tilde{f} \circ g  \:$ is $\:  \mathbf{1}_k\,\mbox{-}\,\!\!\uparrow  \:,$
 too. For $\:  k <  d  \:$ we apply Lemma 2, i.e.
 $$
 \{ \gamma \subseteq [d]  \mid \,|\gamma | \geq k \}   = \!\pbigcup _{\!\!\!j} \langle \sigma_j , \tau _j \rangle
 $$
 is a disjoint union of set intervals, where $\:  \sigma _j \subseteq \tau_j , \sigma _j \sim \tau_j  \:$ and
 $\:  | \sigma_j  | = k  \:$ for each $\:  j  \:.$ Remember that $\: \sigma \sim \tau  \:$ means
 $\:  \max\limits_\sigma x_i = \max\limits_\tau x_i  \:$ (in $\: \R^k  \:).$ Since by (3)
 $$
 \sum_{ \alpha \in \langle \sigma_j, \tau_j \rangle} ( \Delta _\alpha f ) ( \emptyset ) = ( \Delta _{\sigma_j} \circ
 E_{\tau _j \smallsetminus \sigma_j} ) (f) (\emptyset) = ( \Delta _{\sigma_f} f ) ( \tau _j \smallsetminus \sigma_j ) \geq 0
 $$
 (because of  $\:  | \sigma_j | = k ),  \:$ we get
 $$
 \begin{array}{lcl}
   \dis \tilde{f} \circ g &  =  &  \dis \sum_{{\alpha\subseteq [d] \atop | \alpha | < k }} ( \Delta _ \alpha f ) ( \emptyset ) \cdot 1_{[\max\limits_\alpha x_i, \infty] \cap A }    \\[2.5ex]

   & + &  \dis \sum_{ j} ( \Delta _ {\sigma_j} f ) ( \tau_j \smallsetminus \sigma_j ) \cdot 1 _{[\max\limits _{\sigma_j} x_i, \infty ]   \cap A}
   \end{array}
   $$
   which is $\:  \mathbf{1}_k\,\mbox{-}\,\!\!\uparrow \:.$  \\

   In the next step we let $\:  g _1 \:$ be the d.f. of some probability measure with finite support in
   $\:  A  \:,$ say $\:  \sum^n _{\ell = 1} \lambda _\ell \varepsilon _{a_{1, \ell}} \:$ with $\:  \lambda _\ell \geq 0 , \sum _{ \ell} \lambda _\ell = 1 \:$ and
   $\:   a_{1, \ell}  \in A  \:.$ Since $\:  \tilde{f} \:$ is affine as a function of $\:  x_1 \:,$
   $$
   \tilde{f}  \circ g = \sum_\ell \,  \lambda_\ell  \, \tilde{f}  \circ ( g_{1, \ell} , g_2, \ldots , g_d )
   $$
   is again $\: \mathbf{1}_k\,\mbox{-}\,\!\!\uparrow \:.$ This procedure is then repeated for $\:  g_2, g_3, \ldots , g_d  \:,$ showing our result to be true if each $\:  g_i  \:$
   is the d.f. of some probability measure with finite support in $\:  A  \:.$  \\

   Invoking Lemma 3 we may extend the validity to $\:  \mathbf{1}_k\,\mbox{-}\,\!\! \uparrow \:$ functions  $\:  g_1, \ldots , g_d  \:$ for which
   $\: c_i := \sup g_i (A) = 1  \:$ for each $\:  i  \:,$ making use also of the continuity of $\:  \tilde{f} \:.$ \\

   In general we have $\:  c_i \in [0,1] \:,$ where we may assume $\:  c_i > 0  \:$ for each $\:  i  \:.$ Then $\:  \varphi ( x) := \tilde{f} ( c \odot x )  \:$
   is still multilinear and $\:  \mathbf{1}_k\,\mbox{-}\,\!\!\uparrow \:, $ so that
   $$
    \varphi \circ ( g_1 / c_1, \ldots g_d / c_d ) = \tilde{f} \circ g
    $$
   is  $\:  \mathbf{1}_k \,\mbox{-}\,\!\!\uparrow \:,$ thereby finishing our proof. $  \ \Box$  \\

   Theorem 1 deals with fully $k$-increasing functions, generalizing the case $\:  k = 2 \:$ of increasing super-modular
   functions. In \cite{Ch} fully $k$-alternating functions are dealt with, for which we offer the following general result: \\

   \textbf{Theorem 2}. \ \textit{If  in the situation of Theorem} 1 \textit{the functions} $\:  g_1, \ldots , g_d  \:$ \textit{and} $\:  f  \:$
   \textit{are} $\: \mathbf{1}_k\,\mbox{-}\,\!\!\updownarrow \:,$ \textit{then so is}  $\:  h \:.$ \\

   \textbf{Proof}.  \ We make use of the very close direct connection between $\:  \mathbf{n}\,\mbox{-}\,\!\!\uparrow \:$ and
   $\:  \mathbf{n}\,\mbox{-}\,\!\!\updownarrow \:$ functions in full generality --- see \cite{Re3}, Remark (d):
   $$ \begin{array}{l} \dis
   \varphi : A \lra \R \mbox{ is } \mathbf{n}\,\mbox{-}\,\!\!\updownarrow \;  \Longleftrightarrow  \; - \varphi (-\cdot) \mbox{ is }
   \mathbf{n}\,\mbox{-}\,\!\!\uparrow \: \mbox{ on } \:  - A  \\[2ex]

    \Longleftrightarrow c - \varphi ( - \cdot ) \mbox{ is }   \mathbf{n}\,\mbox{-}\,\!\!\uparrow  \mbox{ on }  - A \quad \forall \; c \in \R
    \end{array}
    $$
    where in our situation (i.e. $\:  n_j \in \{ 0, 1 \} \: \forall _j  ) \quad A = \prod^k_{ j = 1} A_j  \:$
    with arbitrary non-empty subsets $\:  A_j \subseteq \overline{\R} \:.$ We apply this to $\:  g_1 , \ldots , g_d  \:$ and to $\:  f  \:: $
    $$ \begin{array}{l} \dis
    g_1, \ldots , g_d : A \lra [0,1] \mbox{ are } \mathbf{1}_k\,\mbox{-}\,\!\!\updownarrow \: \mbox{ and } \: f : \{ 0, 1 \} ^d \lra
    [0,1] \: \mbox{ is } \:  \mathbf{1}_k\,\mbox{-}\,\!\!\updownarrow \\[1.5ex]

     \dis \Longleftrightarrow 1 - g_i (-\cdot) : - A \lra [0,1] \: \mbox{ is } \:  \mathbf{1}_k\,\mbox{-}\,\!\!\uparrow \; \forall \: i  \\[1.5ex]

     \hspace*{1cm} \mbox{ and }  \: 1 - f (-\cdot ) : \{ -1, 0 \} ^d \lra [0,1] \: \mbox{ is } \: \mathbf{1}_k\,\mbox{-}\,\!\!\uparrow \: ,
     \end{array}
     $$
     where the last statement is equivalent with
     $\: 1 - f ( \mathbf{1}_d-\cdot) \:$  being $\: \mathbf{1}_k\,\mbox{-}\,\!\!\uparrow  \:$ on $\: \{ 0, 1 \} ^d  \:.$
     By Theorem 1
     $$
     1- f ( \mathbf{1}_d - (\mathbf{1}_d - g (-\cdot ))) = 1 - f \circ g ( - \cdot ) \: \mbox{ is } \: \mathbf{1}_k\,\mbox{-}\,\!\!\uparrow \: ,
     $$
     or, equivalently, $\:  f \circ  g  \:$ is $\:  \mathbf{1}_k\,\mbox{-}\,\!\!\updownarrow \:. \ \Box $ \\

     \textbf{Remark 1}. \ In \cite{Ch} the functions $\:  g_i  \:$ may be defined on $\: \{ 0, 1 \} ^\ell  \:$ for
     $\:  \ell \geq k  \:.$ This is of course only superficially more general, since by the very definition, being fully
     $k$-increasing or alternating, only $\:  k  \:$ variables are considered simultaneously. \\

     \textbf{Remark 2}. \ For $\:  k = d  \:$ the polynomial $\:  \tilde{f} \:$ has but non-negative coefficients (cf. (7) above), hence not only
     $\:  \tilde{f}  \circ ( g_1, \ldots , g _d) \:$ is $\:  \mathbf{1}_d\,\mbox{-}\,\!\!\uparrow \:,$
     but even $\:  \tilde{f} \circ ( g _1 \times \ldots \times g_d )  \:$ is $\:  \mathbf{1}_{d^2}\,\mbox{-}\,\!\!\uparrow  \:$
     on $\:  A^d  \:.$ For $\:   k < d  \:$ this cannot be expected: take $\:  k = 1, d = 2 , f : \{ 0, 1 \} ^2 \lra \R  \:$
     defined by $\:  f ( 0, 0 ) = 0, f (1,0) = f (0,1) = 1 = f (1,1) \:.$ Then $\:  f  \:$ (and  $ \: \tilde{f}) \:$  is increasing, but
     not $\:  \mathbf{1}_2\,\mbox{-}\,\!\!\uparrow  \:,$ we  have
     $\:  \tilde{f} (x_1, x_2 ) = x_1 + x_2 - x_1 x_2  \:.$ For increasing functions $\:  g_1, g_2  : [0,1]  \lra [0,1] \:$ the composed,
      map $\: \tilde{f} \circ (g_1, g_2 ) = g_1 + g_2 - g_1 g_2  \:$ is still increasing, however the bivariate
      $$
       ( \tilde{f} \circ ( g_1 \times g_2 )) (s,t) = g_1 (s) + g_2 (t) - g_1    (s)  g _2 (t)
       $$
       is not in general $\:  \mathbf{1}_2\,\mbox{-}\,\!\!\uparrow \:,$
       for ex. in the case $\:  g_1 = g_2 = id  \:,$ with $\:  ( \tilde{f} \circ ( g_1 \times g_2 )) (s,t) = s + t - st  \:.$ \\

       \textbf{Example 2}.  \ In Example 1 we considered the case  $\:  k = 2 , d = 3  \:$ and $\:  f (\alpha ) := | \alpha | \vee 1  \:,$
       with
       $$
       \tilde{f} (x) = 1 + x_1 x_2 + x_1 x_3 + x_2 x_3 - x_1 x_2 x_3   \:
       $$
       $ f  \:$ and $\:  \tilde{f} \:$ are  $\:  \mathbf{1}_2\,\mbox{-}\,\!\!\uparrow  \:.$ \\

       If now $\:  g_1, g_2, g_3  \:$ are $\:  \mathbf{1}_2\,\mbox{-}\,\!\!\uparrow \:,$ with values in $\:  [0,1] \:,$
       then by Theorem 1
       $$
       \tilde{f} \circ (g_1, g_2, g_3) = 1 + g_1 g_2 + g_1 g_3 + g_2 g_3 - g_1 g _2 g_3
       $$
       is also $\:  \mathbf{1}_2\,\mbox{-}\,\!\!\uparrow \:,$ hence a distribution function in case $\:  g_1, g_2, g_3  \:$ are of this type. Proving this
       directly, say for differentiable $\:  g_i  \:,$ should be possible, but is certainly cumbersome. \\

       \textbf{Example 3}. \ Again $\:   k = 2 , d = 3 \:.$ Let $\:  f : {\cal P} ([3]) \lra \R  \:$ be given by $\:  f ( \emptyset) = 0,
       f ( \{ i \} ) = 2 \;  ( i = 1, 2, 3), f ( \{ i ,j \} ) = 4 \;  ( i \neq j )  \:$ and $\:  f ( [3]) = 5 \:.$ Then $\:  f  \:$
       is easily seen to be $\:  \mathbf{1}_2\,\mbox{-}\,\!\!\updownarrow \:,$  but not $\:  \mathbf{1}_3\,\mbox{-}\,\!\!\updownarrow\:.$ Here
       $$
       \tilde{f} (x) = 2 ( x_1 + x_2 + x_3 ) -  x_1 x_2 x_3  \: ,
       $$
       and for  increasing submodular $\:  g_1, g_2, g_3  \:$ with values in $\:  [0,1] \:$ also
       $\:  \tilde{f} \circ ( g_1, g_2, g_3 )  \:$ is again increasing and submodular by Theorem 2. \\

\textbf{Remark 3}. \ A natural question is to know which (univariate) functions $\:  \varphi  \:$ \glqq operate\grqq \ on fully $k$-increasing (resp.\ alternating) functions, i.e.\ have
the property that $\:  \varphi \circ f  \:$ is $\:  \mathbf{1}_k\,\mbox{-}\,\!\!\uparrow (\updownarrow) \:$ whenever $\:  f  \:$ is $\:  \mathbf{1}_k\,\mbox{-}\,\!\!\uparrow (\updownarrow) \:,$
supposing of course $\:  \varphi \circ f   \:$ to be defined. The answer is provided by Theorem 12 in \cite{Re2}, later  (in \cite{Re3}, p. 250) called \textit{Monotone Composition Theorem}: if
$\:  \varphi  \:$ is $\:  k\,\mbox{-}\,\!\!\uparrow ( \updownarrow) \:$ and $\:  f  \:$ is $\: \mathbf{1}_k\,\mbox{-}\,\!\! \uparrow (\updownarrow)  \:$ then also
$\:  \varphi \circ f  \:$ is $\: \mathbf{1}_k\,\mbox{-}\,\!\!\uparrow ( \updownarrow)  \:.$ For a pseudo-Boolean function $\:  f  \:$ on $\: \{ 0, 1 \} ^d  \:$ we saw
in Lemma 1 that $\:  f  \:$ is $\: \mathbf{1}_k\,\mbox{-}\,\!\!\uparrow \:$ iff $\:  \tilde{f} \:$ is (on $\: [0,1]^d ), \:$
and this is likewise true for fully $k$-alternating $\:  f \:.$ Hence for $\:  k\,\mbox{-}\,\!\!\uparrow (\updownarrow)  \:\varphi  \:$ also
$\:  \varphi \circ f, \varphi \circ \tilde{f} \:$ and $\:  ( \varphi \circ f ) ^\sim  \:$ are $\:  \mathbf{1}_k\,\mbox{-}\,\!\! \uparrow ( \updownarrow)\:,$
the latter two being different in general. In Example 3 above, with $\:  \varphi (t) := \sqrt{t} \; (k\,\mbox{-}\,\!\!\updownarrow \; \forall \:  k \in \N ) \:,$
we have
$$ \begin{array}{lcl}
\dis ( \sqrt{f})^\sim (x) & = &  \dis  \sqrt{2} ( x_1 + x_2 + x_3) - 2 ( \sqrt{2} - 1 ) (x_1 x_2 + x_1 x_3 + x_2 x_3 ) \\[2ex]

& & \dis + (3 \sqrt{2} + \sqrt{5} - 6  )\, x_1 x_2 x_3 ,
\end{array}
$$
which is $\:  \mathbf{1}_2\,\mbox{-}\,\!\!\updownarrow \:$ on $\:  [0,1]^3 \:,$ as is also
$\: \sqrt{\tilde{f} (x)} \:.$ \\

Note that $\:  \varphi : \R_+ \lra \R_+  \:$ is $\:  2\,\mbox{-}\,\!\!\updownarrow \:$ if $\:  \varphi  \:$
is increasing and concave. And $\: \varphi  \:$ is $\:  k\,\mbox{-}\,\!\!\updownarrow \; \forall \: k \in \N  \:$ iff
it is a socalled Bernstein function. \\

\textbf{Remark 4}. \ When looking at Theorem 1 one might believe that perhaps each $\: \mathbf{1}_k\,\mbox{-}\,\!\!\uparrow \:$ function on
$\: [0,1]^d \:$ has the property shown there for multilinear polynomials $\: \tilde{f} \:$ arising from a
$\: \mathbf{1} _k \,\mbox{-}\,\!\!\uparrow \:$ pseudo-Boolean function $\:  f \:. $ This is not the case:

Let again $\:  k = 2 , d = 3  \:$ and $\:  a := ( \frac{1}{2}, \frac{1}{2}, \frac{1}{2})  \:;$ the d.f. $\: \varphi   \:$ of
$\:  \varepsilon_a  \:$ (on $\: [0,1]^3) \:$ is given by $\:  \varphi = 1 _{[a,1]} \:,$ it is (even) $\:  \mathbf{1}_3\,\mbox{-}\,\!\!\uparrow  \:.$
Let further $\:  g_1 = g_2 = g_3  \:$ be the d.f. of the uniform distribution of $\: [0,1]^2 \:,$  i.e. $\: g_i (s,t) = st \:.$
Then $\:  \varphi \circ ( g_1, g_2, g_3 ) = 1_A  \:$ with
$$
 A := \left\{\{ (s,t) \in [0,1]^2 \mid st \geq \frac{1}{2} \right\} \,
 $$
 and this function is not $\:  \mathbf{1}_2\,\mbox{-}\,\!\!\uparrow \:$ because
 $$
  \left( \Delta ^{ (1,1)}_ { (\frac{1}{2}, \frac{1}{2})} 1 _A \right)  \left( \frac{1}{2}, \frac{1}{2} \right) = - 1  \: .
  $$

\vspace*{.3cm}

\bibliographystyle{plain}

\end{document}